\documentclass[11pt]{amsart}
\usepackage{mathrsfs}
\usepackage{amssymb}

\usepackage{titletoc}
\usepackage{amscd}
\usepackage{amsmath, amssymb}
\usepackage{amsfonts}
\usepackage[colorlinks,linkcolor=blue,citecolor=blue, pdfstartview=FitH]{hyperref}

\numberwithin{equation}{section}

\theoremstyle{plain}
\newtheorem{thm}{Theorem}[section]
\newtheorem{theorem}[thm]{Theorem}
\newtheorem{lemma}[thm]{Lemma}
\newtheorem{corollary}[thm]{Corollary}
\newtheorem{proposition}[thm]{Proposition}

\theoremstyle{definition}

\newtheorem{question}[thm]{Question}

\newtheorem{remark}[thm]{Remark}

\newtheorem{definition}[thm]{Definition}

\newtheorem{example}[thm]{Example}

\newtheorem{defn-thm}[thm]{Definition-Theorem}
\newtheorem{conjecture}[thm]{Conjecture}

\newcommand{\sO}{{\mathcal O}}

\newcommand{\C}{{\mathbb C}}

\renewcommand{\H}{{\mathbb H}}

\renewcommand{\P}{{\mathbb P}}

\newcommand{\R}{{\mathbb R}}

\renewcommand{\S}{{\mathbb S}}

\newcommand{\Z}{{\mathbb Z}}

\newcommand{\qtq}[1]{\quad\mbox{#1}\quad}
\newcommand{\bp}{\bar{\partial}}
\newcommand{\Om}{\Omega}

\newcommand{\ts}{\otimes}

\newcommand{\btheorem}{\begin{theorem}}
\newcommand{\etheorem}{\end{theorem}}
\newcommand{\bproposition}{\begin{proposition}}
\newcommand{\eproposition}{\end{proposition}}
\newcommand{\bdefinition}{\begin{definition}}
\newcommand{\edefinition}{\end{definition}}
\newcommand{\bcorollary}{\begin{corollary}}
\newcommand{\ecorollary}{\end{corollary}}
\newcommand{\bproof}{\begin{proof}}
\newcommand{\eproof}{\end{proof}}
\newcommand{\bremark}{\begin{remark}}
\newcommand{\eremark}{\end{remark}}
\newcommand{\eexample}{\end{example}}
\newcommand{\bexample}{\begin{example}}
\newcommand{\la}{\langle}
\newcommand{\elemma}{\end{lemma}}
\newcommand{\blemma}{\begin{lemma}}
\newcommand{\ra}{\rangle}
\newcommand{\sq}{\sqrt{-1}}

\newcommand{\p}{\partial}

\renewcommand{\bar}{\overline}

\renewcommand{\phi}{\varphi}

\newcommand{\ee}{\end{eqnarray*}}
\newcommand{\be}{\begin{eqnarray*}}

\newcommand{\beq}{\begin{equation}}
\newcommand{\eeq}{\end{equation}}

\newcommand{\bd}{\begin{enumerate}}
\newcommand{\ed}{\end{enumerate}}

\renewcommand{\hat}{\widehat}
\renewcommand{\tilde}{\widetilde}

\renewcommand{\>}{\rightarrow}

\usepackage{fancyhdr}
\begin{document}
\title{Scalar curvature, Kodaira dimension and $\hat A$-genus}

 \makeatletter
\let\uppercasenonmath\@gobble
\let\MakeUppercase\relax
\let\scshape\relax
\makeatother

\author{Xiaokui Yang}
\date{}
\address{{Address of Xiaokui Yang: Morningside Center of Mathematics, Institute of
        Mathematics, Hua Loo-Keng Key Laboratory of Mathematics,
        Academy of Mathematics and Systems Science,
        Chinese Academy of Sciences, Beijing, 100190, China.}}
\email{\href{mailto:xkyang@amss.ac.cn}{{xkyang@amss.ac.cn}}}

\maketitle

\begin{abstract} Let $(X,g)$ be a compact Riemannian manifold with
quasi-positive Riemannian scalar curvature. If there exists a
complex structure $J$ compatible with $g$, then the canonical bundle
$K_X$ is not pseudo-effective and the Kodaira dimension
$\kappa(X,J)=-\infty$. We also introduce the complex Yamabe number
$\lambda_c(X)$ for compact complex manifold $X$, and show that if
$\lambda_c(X)>0$, then $\kappa(X)=-\infty$; moreover, if $X$ is also
spin, then the Hirzebruch $A$-hat genus  $\hat A(X)=0$.
\end{abstract}

\small{\setcounter{tocdepth}{1} \tableofcontents}

\section{Introduction}

This is a continuation of our previous paper \cite{Yang17}, and we
investigate the geometry of Riemannian scalar curvature
on compact complex manifolds. \\

 The existences of various positive scalar curvatures are obstructed. For instance,  it is well-known that, if a compact Hermitian manifold
 has
  quasi-positive \emph{Chern scalar curvature},  then the Kodaira dimension
  is $-\infty$. On the other hand, a classical result of Lichnerowicz
says that if a compact Riemannian spin manifold has quasi-positive
\emph{Riemannian scalar curvature}, then the $\hat A$-genus is zero.
The first main result of our paper is
 \btheorem\label{A} Let
$(X,g)$ be a compact Riemannian manifold with quasi-positive
Riemannian scalar curvature. If there exists a complex structure $J$
compatible with $g$, then the canonical bundle $K_X$ is not
pseudo-effective and $\kappa(X,J)=-\infty$. \etheorem

\noindent Here quasi-positive means  non-negative everywhere and
strictly positive at some point.  As it is well-known, the
positivity of the Riemannian scalar curvature of $(X,J,g)$ can not
imply that of the Chern scalar curvature.  As a border line case, we
obtain:

\btheorem\label{B} Let $(X,g)$ be a compact Riemannian manifold with
zero Riemannian scalar curvature. Suppose there exists a complex
structure $J$ compatible with $g$. If $\kappa(X,J)\geq 0$, then
$(X,\omega_g)$ is a K\"ahler Calabi-Yau manifold and
$\emph{Ric}(\omega_g)=0$.
 \etheorem
\noindent The proofs of Theorem \ref{A} and Theorem \ref{B} rely on
several observations in our previous paper \cite{Yang17} and  a new
scalar curvature relation in  Theorem \ref{T}.

 Note  that, on K\"ahler Calabi-Yau surfaces (e.g. $K3$ surfaces, bi-elliptic
surfaces), there is no Riemannian metrics with quasi-positive
 scalar curvature (e.g. \cite[Theorem~A]{L}). However, by Stolz's solution (\cite[Theorem~A]{S}) to the Gromov-Lawson conjecture,  on a simply
connected K\"ahler Calabi-Yau manifold $X$ with holonomy group
$SU(2m+1)$, there do exist Riemannian metrics with quasi-positive
scalar curvature. On the contrary, as an application of Theorem
\ref{A} and Theorem \ref{B}, we show the Riemannian metrics with
quasi-positive scalar curvature are not compatible with the
Calabi-Yau complex structures, and more generally

\bcorollary\label{C} On a compact complex Calabi-Yau manifold $X$
with torsion canonical bundle $K_X$, there is no Hermitian metric
with quasi-positive Riemannian scalar curvature. Moreover, if $X$ is
also non-K\"ahler, then there is no Hermitian metric with
non-negative Riemannian scalar curvature.

\ecorollary

\noindent It is well-known that all compact K\"ahler Calabi-Yau
manifolds have torsion canonical bundle. On the other hand, many
non-K\"ahler Calabi-Yau manifolds also have torsion canonical
bundle. For instance, the connected sum $\#_k(\S^3\times \S^3)$ with
$k\geq 2$
(\cite{LT}).\\

On a compact complex manifold $X$ of complex dimension $n\geq 2$, we
introduce the complex Yamabe number $\lambda_c(X)$: \beq
\lambda_c(X)=\sup_{ g\text{ is Hermitian}}\ \ \ \ \ \inf_{\tilde g\
\text{is conformal to}\ g}\frac{\int_X s_{\tilde g} dV_{\tilde
g}}{\left(\int_X dV_{\tilde {g}}\right)^{1-\frac{1}{n}}},\label{ya}
\eeq where $s_{\tilde g}$ is the Riemannian scalar curvature of
$\tilde g$. Note that in (\ref{ya}), if the supremum is taken over
all Riemannian metrics, then it is the classical Yamabe number
$\lambda(X)$ in conformal geometry. Our second main result is

 \btheorem\label{E} Let $X$ be a compact
complex manifold. If $\lambda_c(X)>0$, then $\kappa(X)=-\infty$.
Moreover, if $X$ is also spin, then $\hat A(X)=0$. \etheorem
\noindent According to the results of Gromov-Lawson \cite{GL} and
Stolz \cite{S},  on a simply connected K\"ahler Calabi-Yau manifold
$X$ with $\dim_\C X\geq 3$, one
has $\lambda(X)>0$ and $\hat A(X)=0$. However, we have $\lambda_c(X)\leq 0$ by Theorem \ref{E}.\\

As motivated by Theorem \ref{A}, Theorem \ref{B}, Theorem \ref{E},
various conjectures described  in \cite[Section~4]{Yang17} and
classical works by Schoen-Yau \cite{SY1, SY2, SY3}, Gromov-Lawson
\cite{GL}, Stolz \cite{S} and LeBrun \cite{L}  (see also Zhang
\cite{Z}), we propose

\begin{conjecture}\label{con} Let $X$ be a compact K\"ahler manifold with
$\kappa(X)=-\infty$. If $X$ has a spin structure, then $\hat
A(X)=0$.
\end{conjecture}
 \noindent Note that Conjecture \ref{con} holds when
$\dim_\C X=2$ or $2m+1$.\\

 Finally, we would like to describe some  applications of Theorem
 \ref{A}.

\bproposition\label{D} Let $X$ be a compact K\"ahler threefold. If
there exists a Hermitian metric with quasi-positive Riemannian
scalar curvature, then $X$ is uniruled, i.e. $X$ is covered by
rational curves. \eproposition

\noindent According to the uniruledness conjecture (e.g.
\cite[Conjecture]{BDPP13}), Proposition \ref{D} should be true
on higher dimensional compact K\"ahler manifolds.\\

It is a long-standing open problem to determine whether the
six-sphere $\S^6$ admits  a complex structure or
not\footnote{Recently, it was announced in \cite{A}  by Sir Michael
Atiyah that there is no complex structure on $\S^6$.}. Now assuming
$X:=\S^6$ has a complex structure $J$. As  pointed out in
\cite[p.~122]{HKP}, it is not at all clear whether
$\kappa(X,J)=-\infty$, and proving this would seem to be as
complicated as to show that there are no divisors on $X$ at all. It
is obvious that $c_1(X)=0\in H^2(X,\Z)$ and it is also proved in
\cite{T} that $c_1^{\text{BC}}(X,J)\neq 0$ and in particular, $K_X$
is not holomorphically torsion. For more related discussions, we
refer to \cite{An}. Let $\mathscr S$ be the space of Riemannian
metrics with non-negative scalar curvature. We have

\btheorem\label{S6} If there exists a complex structure $J$ which is
compatible with some $g\in \mathscr S$, then $K_X$ is not
pseudo-effective and
$$\kappa(X,J)=-\infty.$$
 \etheorem

\noindent It is known that there is no complex structure compatible
with metrics in a small neighborhood of the round metric on $\S^6$
(e.g. \cite{L2,LY14,Tang, BHL99}).\\

\noindent\textbf{Acknowledgement.} The author would like to thank
Bing-Long Chen, Fei Han,
 Tian-Jun Li,
 Ke-Feng Liu, Jian-Qing Yu, Bai-Lin Song, Song Sun,  Valentino Tosatti, Shing-Tung Yau and Fang-Yang Zheng for
valuable comments and discussions. This work was partially supported
by China's Recruitment
 Program of Global Experts, National Center for Mathematics and Interdisciplinary Sciences,
 Chinese Academy of Sciences.

\vskip 2\baselineskip

\section{Preliminaries}

\subsection{Ricci curvature on complex manifolds}
 Let $(X,\omega_g)$ be a compact Hermitian manifold. Locally, we write
$$\omega_g=\sq g_{i\bar j}dz^i\wedge d\bar z^j$$
  The (first Chern-)Ricci form $\text{Ric}(\omega_g)$ of
$(X,\omega_g)$ has components
$$R_{i\bar j}=-\frac{\p^2\log\det(g_{k\bar \ell})}{\p z^i\p\bar z^j}$$
which also represents the first Chern class $c_1(X)$ of the complex
manifold $X$ (up to a constant). The Chern scalar curvature
$s_\text{C}$ of $(X,\omega_g)$ is given by \beq
s_\text{C}=\text{tr}_{\omega_g}\text{Ric}(\omega_g)=g^{i\bar j}
R_{i\bar j}. \eeq The total Chern scalar curvature of $\omega_g$ is
\beq \int_X s_\text{C}\cdot  \omega_g^n=n
\int\text{Ric}(\omega_g)\wedge \omega_g^{n-1},\eeq where $n$ is the
complex dimension of $X$.

\subsection{The Bott-Chern classes}
The Bott-Chern
  cohomology  on a compact complex manifold
  $X$ is given by
  $$ H^{p,q}_{\mathrm{BC}}(X):= \frac{\text{Ker} d \cap \Om^{p,q}(X)}{\text{Im} \p\bp \cap \Om^{p,q}(X)}.$$

\noindent Let $\mathrm{Pic}(X)$ be the set of holomorphic line
bundles over $X$. As similar as the first Chern class map
$c_1:\mathrm{Pic}(X)\>H^{1,1}_{\bp}(X)$, there is a \emph{first
Bott-Chern class} map \beq c_1^{\mathrm{BC}}:\mathrm{Pic}(X)\>
H^{1,1}_{\mathrm{BC}}(X).\eeq Given any holomorphic line bundle
$L\to X$ and any Hermitian metric $h$ on $L$, its curvature form
$\Theta_h$ is locally given by $-\sq\p\bp\log h$.  We define
$c_1^{\mathrm{BC}}(L)$ to be the class of $\Theta_h$ in
$H^{1,1}_{\mathrm{BC}}(X)$.  For a complex manifold $X$,
$c_1^{\mathrm{BC}}(X)$ is defined to be
$c_1^{\mathrm{BC}}(K^{-1}_X)$ where $K_X^{-1}$ is the anti-canonical
line bundle.

\subsection{Special manifolds} Let $X$ be a compact complex
manifold.

 \bd \item A Hermitian metric $\omega_g$ is called a Gauduchon metric if
$\p\bp\omega_g^{n-1}=0$. It is proved by Gauduchon (\cite{Ga2})
that, in the conformal class of each Hermitian metric, there exists
a unique Gauduchon metric (up to  scaling).

\item A Hermitian metric $\omega_g$ is called a K\"ahler metric if
$d\omega_g=0$.

\item $X$ is called a Calabi-Yau manifold if $c_1(X)=0\in
H^2(X,\Z)$.

\ed

\subsection{Kodaira dimension of compact complex manifolds.} The Kodaira dimension $\kappa(L)$ of a line bundle $L$ over
a compact complex manifold $X$ is defined to be
$$\kappa(L):=\limsup_{m\>+\infty} \frac{\log \dim_\C
H^0(X,L^{\ts m})}{\log m}$$ and the \emph{Kodaira dimension}
$\kappa(X)$ of $X$ is defined as $ \kappa(X):=\kappa(K_X)$ where the
logarithm of zero is defined to be $-\infty$. In particular, if
$$\dim_\C H^0(X,K_X^{\ts m})=0$$
for every $m\geq 1$, then $\kappa(X)=-\infty$.

\subsection{Spin manifold and $\hat A$-genus.} Let $X$ be a compact
oriented  Riemannian manifold. It is called a spin manifold, if it
admits a spin structure, i.e. its second Stiefel-Whitney class
$w_2(X)=0.$ It is well-known that all compact Calabi-Yau manifolds
are spin.

The definition of the $\hat A$-genus of a compact oriented
Riemannian manifold $X$ is as follows. Let $\hat A_i(p_1,\cdots,
p_i)$ be the multiplicative sequence of polynomials in the
Pontryagin classes $p_i$ of $X$ belonging to the power series
$$\frac{\frac{1}{2}\sqrt{z}}{\sinh\left(\frac{1}{2}\sqrt{z}\right)}=1-\frac{1}{24}z+\frac{7}{2^7\cdot 3^2\cdot 5}z^2+\cdots.$$
The first few terms are
$$\hat A_1(p_1)=-\frac{1}{24}p_1,\ \ \ \hat A_2(p_1,p_2)=\frac{1}{2^7\cdot 3^2\cdot 5}\left(-4p_2+7p_1^2\right).$$
The $\hat A$-genus, $\hat A(X)$ is by definition the real number
$(\sum_i \hat A_i(p_1,\cdots, p_i))[X]$, where $[X]$ means
evaluation of the cohomology class on the fundamental cycle of $X$.
Since $p_i\in H^{4i}(X,\Z)$, $\hat A(X)$ is zero unless $\dim_\R
X\equiv 0 (\text{mod}\ 4)$. Moreover, if $X$ is a spin manifold,
$\hat A(X)$ is an integer.  The following result is well-known (for
more historical explanations, we refer to \cite{LM, Wu} and the
reference therein):

\blemma On  a compact spin manifold $X$, if it admits a Riemannian
metric with quasi-positive scalar curvature, then $\hat A(X)=0$.
\elemma

\noindent For more necessary background materials, we refer to
\cite{LM,LY12, LY14, Yang16, Yang17, LY17} and the references
therein. \vskip 2\baselineskip
\section{The  Riemannian scalar curvature  and Kodaira dimension}

Let $(X,\omega)$ be a compact Hermitian manifold. We first give
several computational results.

\blemma For any smooth real valued function $f\in C^\infty(X,\R)$,
we have \beq \bp^*(f\omega)=f\bp^*\omega+\sq\p f.\label{c1}\eeq
\bproof For any smooth $(1,0)$-form $\eta\in \Gamma(X,T^{*1,0}X)$,
we have the global inner product \be
\left(\bp^*(f\omega),\eta\right)&=&\left(f\omega,\bp\eta\right)
=\left(\omega,
f\bp\eta\right)\\
&=&\left(\omega,\bp(f\eta)\right)-\left(\omega,\bp f\wedge
\eta\right)\\
&=&\left(f\bp^*\omega,\eta\right)-\left(\omega,\bp f\wedge
\eta\right)\\
&=&\left(f\bp^*\omega,\eta\right)+\sq\left(\p f, \eta\right) \ee
where the last identity follows from the fact that $f$ is real
valued. \eproof \elemma

\blemma For any $(1,0)$ form $\eta$ and real valued function $f\in
C^\infty(X,\R)$, we have \beq \p^*(f\eta)=f\p^*\eta-\la \eta,\p
f\ra.\label{c5}\eeq \bproof For any smooth function $\phi\in
C^\infty(X)$, we have  \be
\left(\p^*(f\eta),\phi\right)&=&\left(f\eta,\p\phi\right)=\left(\eta,
f\p\phi\right)=\left(\eta, \p(f\phi)-\p f\cdot \phi\right)\\
&=&\left(f\p^*\eta,\phi\right)-\left(\la\eta, \p f\ra,
\phi\right)\ee and we obtain (\ref{c5}). \eproof \elemma

Let $\omega_f=e^f\omega$ for some $f\in C^\infty(X,\R)$. We denote
by $\bp_f^*, \bp^*_f$ and $\p^*,\bp^*$ the adjoint operators taking
with respect to $\omega_f$ and $\omega$ respectively. The local and
global inner products with respect to $\omega$ and $\omega_f$ are
indicated by $\la\bullet,\bullet\ra, (\bullet, \bullet)$ and
$\la\bullet,\bullet\ra_f, (\bullet, \bullet)_f$ respectively.

\blemma  For any $(1,0)$ form $\eta$ and real valued function $f\in
C^\infty(X,\R)$, we have \beq \p_f^*\eta
=e^{-f}\left[\p^*\eta-(n-1)\left\la \eta,\p
f\right\ra\right].\label{c7}\eeq  \elemma

\bproof For any $\phi\in C^\infty(X)$, we have \be \left(\p_f^*\eta,
\phi\right)_f=\left(\eta,\p \phi\right)_f=(e^{(n-1)f}\eta,
\p\phi)=\left(\p^*\left(e^{(n-1) f}\eta\right),\phi\right)\ee where
the second identity holds since $\eta$ is a $(1,0)$-form. By
(\ref{c5}), we obtain
$$\left(\p_f^*\eta,
\phi\right)_f=\left(e^{(n-1)f}\p^*\eta,\phi\right)-\left(\la\eta, \p
e^{(n-1) f}\ra, \phi\right).$$ Hence,
$$\left(\p_f^*\eta,
\phi\right)_f=\left(e^{-f}\p^*\eta,\phi\right)_{f}-\left((n-1)e^{-f}\la\eta,
\p f\ra, \phi\right)_f,$$ which verifies (\ref{c7}). \eproof

\blemma We have \beq \bp^*_f\omega_f=\bp^*\omega+(n-1)\sq\p
f.\label{c4}\eeq \bproof For any $\eta\in\Gamma(X,T^{*1,0}X)$, we
have \be
\left(\bp_f^*\omega_f,\eta\right)_{f}&=&\left(\omega_f,\bp\eta\right)_{f}=\left(e^{(n-1)f}\cdot\omega, \bp\eta\right)\\
&=&\left(\bp^*\left(e^{(n-1)f}\cdot\omega\right),\eta\right).\ee Now
by (\ref{c1}), we have \be
\left(\bp_f^*\omega_f,\eta\right)_{f}&=&\left(e^{(n-1)f}\left[
\bp^*\omega+(n-1)\sq\p f\right], \eta\right)\\
&=&\left( \bp^*\omega+(n-1)\sq\p f, \eta\right)_{f}\ee since $\eta$
is a $(1,0)$ form. Therefore, we obtain (\ref{c4}).
 \eproof
 \elemma

\blemma We have \beq
\sq\p_f^*\bp^*_f\omega_f=e^{-f}\left(\sq\p^*\bp^*\omega-(n-1)\left(\Delta_df+\mathrm{tr}_\omega\sq
\p\bp f\right)+(n-1)^2|\p f|^2\right).\label{c6}\eeq \bproof By
formulas (\ref{c5}) and (\ref{c4}), we have \be
\sq\p_f^*\bp_f^*\omega_f&=&\sq \p_f^*\left(\bp^* \omega+(n-1)\sq\p
f\right)\\
&=&e^{-f}\left(\sq\p^*\bp^*\omega-\sq(n-1)\la\bp^*\omega,\p
f\ra-(n-1)\p^*\p f+(n-1)^2|\p f|^2 \right).\ee We also observe that
\beq  \sq\la\bp^*\omega,\p f\ra=\bp^*\bp f+\text{tr}_\omega\sq \p\bp
f.\label{c8}\eeq Indeed, for any text function $\phi\in
C^\infty(X)$, we have \be \left(\sq\la\bp^*\omega,\p
f\ra,\phi\right)&=&\sq\left(\bp^*\omega,\phi\p f\right)\\
&=&\sq\left(\omega,\bp\phi\wedge \p f+\phi\bp\p f\right)\\
&=&(\bp f, \bp\phi)+(\omega, \phi\cdot \sq\p\bp f)\\
&=&(\bp^*\bp f, \phi)+(\text{tr}_\omega\sq \p\bp f, \phi)\ee which
gives (\ref{c8}). Since $\Delta_d f=d^*d f=\bp^*\bp f+\p^*\p f$, we
obtain (\ref{c6}).
 \eproof
 \elemma

\noindent The following observation is one of the key ingredients in
the curvature computations. \blemma Let $(X,\omega)$ be a compact
Hermitian manifold. Then \beq \la
\bp\bp^*\omega,\omega\ra=|\bp^*\omega|^2- \sq
\p^*\bp^*\omega.\label{c2} \eeq In particular, if $\omega$ is a
Gauduchon metric, we have \beq \la
\bp\bp^*\omega,\omega\ra=|\bp^*\omega|^2.\label{c3}\eeq

 \bproof
For any smooth real valued function $\phi\in C^\infty(X,\R)$, we
have \be\left(\la \bp\bp^*\omega,\omega\ra,
\phi\right)&=&\left(\bp\bp^*\omega,\phi\omega\right)=\left(\bp^*\omega,
\bp^*(\phi\omega)\right)\\
&=&\left(\bp^*\omega, \phi\bp^*\omega+\sq\p
\phi\right)\\&=&\left(|\bp^*\omega|^2,
\phi\right)+\left(-\sq\p^*\bp^*\omega, \phi\right) \ee where we use
formula (\ref{c1}) in the second identity. Since $\phi$ is an
arbitrary smooth real function, we obtain (\ref{c2}). If $\omega$ is
Gauduchon, i.e. $\p\bp\omega^{n-1}=0$, we have $\p^*\bp^*\omega=0$,
and so (\ref{c3}) follows from (\ref{c2}). \eproof \elemma

\bcorollary\label{sca} On a compact Hermitian manifold $(X,\omega)$,
the Riemannian scalar curvature $s$ and the Chern scalar curvature
$s_{\emph{C}}$ are related by \beq
s=2s_{\emph{C}}-2\sq\p^*\bp^*\omega-\frac{1}{2}|T|^2.\label{key}\eeq
where $T$ is the torsion tensor of the Hermitian metric $\omega$.
 \ecorollary
\bproof By Lemma \ref{sc} in the Appendix, we have
$$s=2s_{\text{C}}+\left(\la\bp\bp^*\omega+\p\p^*\omega,\omega\ra-2|\bp^*\omega|^2\right)-\frac{1}{2}|T|^2.$$
 Hence, by formula (\ref{c2}) we obtain (\ref{key}).\eproof

 Let $\omega_f=e^f\omega$ be a smooth Gauduchon
metric (i.e. $\p\bp\omega_f^{n-1}=0$) in the conformal class of
$\omega$ for some  smooth function $f\in C^\infty(X,\R)$.
\btheorem\label{T} The total Chern scalar curvature of the Gauduchon
metric $\omega_f$ is \beq n\int_X \emph{Ric}(\omega_f)\wedge
\omega_f^{n-1}=\int_X e^{(n-1)f} \cdot \left(\frac{s}{2}+\frac{|T|^2
}{4}\right) \omega^{n}+(n-1)^2\|\p f\|^2_{\omega_f}.\label{key2}\eeq
\etheorem

\bproof Indeed, since $\omega_f$ satisfies $\p\bp\omega_f^{n-1}=0$,
we have  \be n\int_X \text{Ric}(\omega_f)\wedge \omega_f^{n-1}&=&
n\int_X
\text{Ric}(\omega)\wedge \omega_f^{n-1}\\
&=&n\int_X e^{(n-1)f}\cdot \text{Ric}(\omega)\wedge
\omega^{n-1}=\int_X e^{(n-1)f} \cdot s_{\text{C}}\cdot
\omega^{n}\\
&=&\int_X e^{(n-1)f} \left(\frac{s}{2}+\frac{|T|^2 }{4}\right)
\omega^{n}+\int_X e^{(n-1)f} \cdot \sq \p^*\bp^*\omega \cdot
\omega^{n}, \ee where we use the scalar curvature relation
(\ref{key}) in the third identity. Since $\omega_f$ is Gauduchon,
and we have $\p_f^*\bp_f^*\omega_f=0$. By formula (\ref{c6}), we
have \beq
\sq\p^*\bp^*\omega=(n-1)\left(\Delta_df+\mathrm{tr}_\omega\sq \p\bp
f\right)-(n-1)^2|\p f|^2.\eeq It is easy to show that
$$\int_X e^{(n-1)f} \mathrm{tr}_\omega\sq \p\bp
f\cdot \omega^n=n\int_X \sq\p\bp f\wedge \omega^{n-1}_f=0$$ and
$$\int_X e^{(n-1)f}|\p f|^2 \omega^n=\|\p f\|^2_{\omega_f}.$$
Moreover, \be \int_X e^{(n-1)f}\Delta_d f \omega^n&=&\left(d^*df,
e^{(n-1)f}\right)\\&=&(n-1)\left(df, e^{(n-1)f}df\right)\\
&=&(n-1)\left(df,df\right)_f\ee since $df$ is a $1$-form. Finally,
we obtain
$$\int_X e^{(n-1)f} \cdot \sq \p^*\bp^*\omega \cdot
\omega^{n}=(n-1)^2\|df\|^2_{\omega_f}-(n-1)^2\|\p
f\|^2_{\omega_f}=(n-1)^2\|\p f\|^2_{\omega_f}.$$ Putting  all
together, we get (\ref{key2}).
 \eproof

\noindent\emph{The proof of Theorem \ref{A}}. Let $\omega$ be the
Hermitian metric of $(g,J)$. Let $\omega_f=e^f\omega$ be a smooth
Gauduchon metric
 in the conformal class of $\omega$. If the Riemannian scalar curvature $s$ of $\omega$ is quasi-positive, then by
formula (\ref{key2}), the total Chern scalar curvature of the
Gauduchon metric $\omega_f$ is strictly positive, i.e. $$n\int_X
\mathrm{Ric}(\omega_f)\wedge \omega_f^{n-1}>0.$$ By \cite[Theorem
~1.1]{Yang17} and \cite[Corollary~3.2]{Yang17}, $K_X$ is not
pseudo-effective and $\kappa(X,J)=-\infty$. \qed

\vskip 2\baselineskip

\noindent\emph{The proof of Theorem \ref{B}}. Suppose $\omega$ is
not a K\"ahler metric, i.e. the torsion $|T|^2$ is not identically
zero,
 then by formula (\ref{key2}), there exists a Gauduchon metric with positive total Chern scalar curvature. Hence, by \cite[Corollary~3.2]{Yang17} we have  $\kappa(X)=-\infty$ which is a contradiction.
 Therefore, $\omega$ is K\"ahler and so in formula (\ref{key2}), $f$ is a constant and $T=0$.
 That means,  $\omega$ is a K\"ahler metric with zero scalar
 curvature. By \cite[Corollary~1.6]{Yang17}, $X$ is a Calabi-Yau
 manifold since $\kappa(X)\geq 0$. By the Calabi-Yau theorem (\cite{Yau78}), there exists a
 K\"ahler Ricci-flat metric $\omega_{\mathrm{CY}}$, i.e. $\mathrm{Ric}(\omega_{\mathrm{CY}})=0$.
 Hence,
 $$\mathrm{Ric}(\omega)=\mathrm{Ric}(\omega)-\mathrm{Ric}(\omega_{\mathrm{CY}})=\sq\p\bp F$$
 where $F=\log\left(\frac{\omega^n_{\mathrm{CY}}}{\omega^n}\right)$.
Since $\omega$ has zero scalar curvature, we have
$$\Delta_\omega F=\mathrm{tr}_\omega\sq\p\bp F=0,$$
which implies $F=\mathrm{const}$ and $\mathrm{Ric}(\omega)=0$.\qed

\bcorollary Let $(X,g)$ be a compact Riemannian manifold with
nonnegative Riemannian scalar curvature. If there exists a complex
structure $J$ which is compatible with $g$, then either \bd\item
$\kappa(X,J)=-\infty$; or
\item $\kappa(X,J)=0$ and $(X,J,g)$ is a K\"ahler Calabi-Yau.
\ed

\bproof If the Riemannian scalar curvature is not identically zero,
or $s$ is identically zero but the metric is not K\"ahler, then the
total Chern scalar curvature in formula (\ref{key2}) is strictly
positive. Hence $\kappa(X,J)=-\infty$. The last possibility is that
$(g,J)$ is a K\"ahler metric with zero scalar curvature. In this
case, it follows from \cite[Corollary~1.6]{Yang17} that either
$\kappa(X)=-\infty$, or  $\kappa(X,J)=0$ and $(X,J,g)$ is a K\"ahler
Calabi-Yau. \eproof

\ecorollary

\bproposition Suppose $X$ is a compact complex manifold with
$c_1^{\mathrm{BC}}(X)\leq 0$.  Then \bd\item there exists a
Hermitian metric with non-positive Riemannian scalar curvature;
\item there is no Hermitian metric with quasi-positive Riemannian scalar
curvature.\ed Moreover, $X$ admits a Hermitian metric $g$ with zero
Riemannian scalar curvature if and only if $(X,g)$ is a K\"ahler
Calabi-Yau.
 \eproposition

\bproof Note that by definition there exists a $d$-closed
non-positive $(1,1)$ form $\eta$ which represents
$c_1^{\mathrm{BC}}(X)$. By \cite[Theorem~1.3]{STW}, there exists a
non-K\"ahler Gauduchon metric $\omega_G$ such that
$$\text{Ric}(\omega_G)=\eta\leq 0.$$
Hence, for any Gauduchon metric $\omega$, \beq \int_X
\text{Ric}(\omega)\wedge \omega^{n-1}=\int_X
\text{Ric}(\omega_G)\wedge \omega^{n-1}\leq 0. \label{key3} \eeq
\noindent (1). Since $\omega_G$ is Gauduchon, by formula
(\ref{key}), we have
$$s=2s_{\mathrm{C}}-\frac{1}{2}|T|^2=2\text{tr}_{\omega_G}\text{Ric}(\omega_G)-\frac{1}{2}|T|^2\leq 0.$$
(2). If there exists a Hermitian metric with quasi-positive
Riemannian scalar curvature, then it induces a Gauduchon metric with
positive total Chern scalar curvature which is a contradiction.\\

 Suppose $X$ admits a Hermitian metric $g$ with zero
Riemannian scalar curvature, then by formulas (\ref{key2}) and
(\ref{key3}), we have $T=0$ and $f=0$, i.e. $(X,\omega_g)$ is a
K\"ahler manifold with zero scalar curvature. Since
$c_1^{\text{BC}}(X)=c_1(X)\leq 0$, we have $\text{Ric}(\omega_g)=0$.
 \eproof

\vskip 1\baselineskip

\noindent \emph{The proof of Corollary \ref{C}}. If $K_X$ is a
torsion, i.e. $K_X^{\ts m}=\sO_X$ for some $m\geq 1$, then
$\kappa(X)=0$. The first part of Corollary \ref{C} follows from
Theorem \ref{A}, and the second part follows from Theorem \ref{B}.
\qed

\vskip 1\baselineskip

\noindent \emph{The proof of Proposition \ref{D}}. By Theorem
\ref{A}, $K_X$ is not pseudo-effective and $\kappa(X)=-\infty$.
Hence by \cite[Corollary~1.2]{Br} or \cite[Corollary~1.4]{HP}, we
conclude $X$ is uniruled. \qed

\vskip 1\baselineskip

\noindent \emph{The proof of Theorem \ref{S6}}. Since $H^2(X,\R)=0$,
the Hermitian metric $(g,J)$ is not K\"ahler.  Then Theorem \ref{S6}
follows from Theorem \ref{A} and Theorem \ref{B}. \qed

\vskip 2\baselineskip

\section{The Yamabe number, $\hat A$-genus and Kodaira dimension}
Let $(X,g_0)$ be a compact Riemannian manifold with real dimension
$2n$. The \emph{Yamabe invariant}  $\lambda(X,g_0)$ of the conformal
class $[g_0]$ is defined as \beq \lambda(X,g_0)=\inf_{g=e^f g_0,
f\in C^\infty(X,\R)}\frac{\int_X s_g dV_g}{\left(\int_X
dV_g\right)^{1-\frac{1}{n}}}\eeq where $s_g$ is the Riemannian
scalar curvature of $g$. Moreover, one can define the Yamabe number
\beq \lambda(X)=\sup_{\text{all Riemannian metric}\  g}
\lambda(X,g).\label{yamabe} \eeq  \noindent As analogous to
(\ref{yamabe}), on a compact complex manifold $X$, one can define
the complex version \beq \lambda_c(X)=\sup_{\text{all Hermitian
metric}\  g} \lambda(X,g).\label{hyamabe} \eeq \btheorem\label{Y}
Let $X$ be a compact complex manifold. If $\lambda_c(X)>0$, then
$\kappa(X)=-\infty$. Moreover, if $X$ is also spin, then $\hat
A(X)=0$. \etheorem \bproof Suppose $\lambda_c(X)>0$, then there
exists a Hermitian metric $g_0$ such that
$$\lambda(X,g_0)=\inf_{g\in [g_0]}\frac{\int_X s_g
dV_g}{\left(\int_X dV_g\right)^{1-\frac{1}{n}}}>0.$$ Let
$\omega_f=e^f \omega_{g_0}$ be a Gauduchon metric in the conformal
class of $\omega_{g_0}$. Hence, $\omega_f$ has positive total
Riemannian scalar curvature
$$\int_X s_f \cdot
\omega_f^n>0.$$
 Moreover, by formula (\ref{key}), the total
Chern scalar curvature of $\omega_f$ is  \beq \int_X
(s_{\mathrm{C}})_f \cdot \omega_f^n =\int_X \frac{s_f}{2} \cdot
\omega_f^n+\frac{1}{4}\int_X |T_f|^2_f\cdot \omega_f^n>0,\eeq where
we use the fact that $\omega_f$ is Gauduchon, i.e.
$\p^*_f\bp^*_f\omega_f=0$. Therefore, the Gauduchon metric
$\omega_f$ has positive total Chern scalar curvature, and by
\cite[Corollary~3.2]{Yang17}, $\kappa(X)=-\infty$. We also have
$\lambda(X)\geq \lambda_c(X)>0$. It is well-known that
$\lambda(X)>0$ if and only if  there exists a Riemannian metric with
positive Riemannian scalar curvature. If $X$ is spin, then $\hat
A(X)=0$ by Lichnerowicz's result. \eproof

Note that on a simply connected K\"ahler Calabi-Yau manifold $X$
with $\dim_\C X=2m+1$, one has $\lambda(X)>0$ and $\hat A(X)=0$.
However,
 $\lambda_c(X)\leq 0$.

\begin{question} On a compact K\"ahler (or complex) manifold $X$,  find sufficient and necessary conditions such that
$\lambda(X)$ and $\lambda_c(X)$ have the same sign, or
$\lambda(X)=\lambda_c(X)$.
\end{question}

\noindent A result along this line is \bcorollary Let $X$ be a
simply compact complex manifold with $\dim_\C X\geq 3$. If
$\lambda_c(X)$ has the same sign as $\lambda(X)$, then
$\kappa(X)=-\infty$. \ecorollary

\bproof By Gromov-Lawson \cite{GL}  and Stolz \cite{S}, if $X$ is a
simply connected complex manifold with $\dim_\C X\geq 3$, then $X$
has a Riemannian metric with positive scalar curvature, hence
$\lambda(X)>0$ and so $\lambda_c(X)>0$. By Theorem \ref{Y}, we
obtain $\kappa(X)=-\infty$. \eproof

\noindent Finally, we want to present a nice result of LeBrun, which
answers Conjecture \ref{con} affirmatively when $X$ is a compact
spin K\"ahler surface (for related works, see also \cite{GLe} and
\cite{L2}):

\btheorem\cite[Theorem~A]{L} Let $X$ be a compact K\"ahler surface,
then \beq
\begin{cases}\lambda(X) > 0\ \ \text{if and only if }\ \
\kappa(X)=-\infty;\\
\lambda(X)=0 \ \ \text{if and only if }\ \ \kappa(X)=0 \qtq{or} 1;\\
\lambda(X)<0\ \ \text{if and only if }\ \ \kappa(X)=2.
\end{cases}\eeq \etheorem

\noindent According to Theorem \ref{A},  Theorem \ref{B}, Theorem
\ref{E} and \cite[Theorem~1.1]{Yang17}, there should be some
analogous results for $\lambda_c(X)$ on compact K\"ahler manifolds,
which will be addressed in future studies. For some related
settings, we refer to \cite{AD, DS, ACS, LU} and the references
therein.

\vskip 2\baselineskip

\section{Examples on compact non-K\"ahler Calabi-Yau surfaces}

In this section, we discuss two special Calabi-Yau surfaces of class
$\text{VII}$. One is the diagonal Hopf surface $\S^1\times \S^3$ and
the other one is the Inoue surface. It is well-known,  they are
non-K\"ahler Calabi-Yau surfaces with Kodaira dimension $-\infty$,
$b_1(X)=1$ and $b_2(X)=0$. A
 straightforward computation shows that   on $\S^1\times
 \S^3$, there are Hermitian metrics with positive (resp. negative, zero)
 Riemannian scalar curvature (\cite[Section~6]{LY14}). We show by the following example
 that the converses of Theorem \ref{A}  and Theorem \ref{E} are   not valid in general:

 \bproposition On  an Inoue surface $X$, it has $\kappa(X)=-\infty$ and $\hat A(X)=0$. However, it  can not support a Hermitian metric with
 non-negative Riemannian scalar curvature. In particular, $\lambda_c(X)\leq
 0$.

 \eproposition

\bproof On each Inoue surface, there exists a smooth Gauduchon
metric with non-positive Ricci curvature. Indeed, let $(w,z)\in
\H\times \C$ be the holomorphic coordinates, then by the precise
definition of each Inoue surface (\cite{DPS, Tel, FTWZ, LY17}), we
know the form
$$\sigma=\frac{dw\wedge dz}{\text{Im}(w)}$$
descends to a smooth  nowhere vanishing $(2,0)$ form on $X$, i.e.
$\sigma\in \Gamma(X,K_X)$. Then it induces a smooth Hermitian metric
$h$ on $K_X$ given by $h(\sigma,\sigma)=1$.
 In the holomorphic frame $e=dw\wedge dz$ of $K_X$, we have
 $$h=h(e,e)=[\text{Im}(w)]^2.$$
It also induces a Hermitian metric $h^{-1}$ on $K_X^{-1}$, and the
curvature of $h^{-1}$ is
$$-\sq \p\bp\log h^{-1}=\sq\p\bp\log [\text{Im}(w)]^2=-\frac{\sq}{2}\frac{dw\wedge d\bar w}{[\text{Im}(w)]^2},$$
which also represents $c^{\mathrm{BC}}_1(X)$. By Theorem
\cite[Theorem~1.3]{STW}, there exists a Gauduchon metric $\omega_G$
with non-positive Ricci curvature
$$\text{Ric}(\omega_G)=-\frac{\sq}{2}\frac{dw\wedge d\bar
w}{[\text{Im}(w)]^2}\leq 0.$$ (Note also that the Riemannian scalar
curvature of $\omega_G$ is strictly negative according to
(\ref{key}).) Hence, for any Gauduchon metric $\omega$, the total
Chern scalar curvature
$$2\int_X\text{Ric}(\omega)\wedge \omega=2\int_X\text{Ric}(\omega_G)\wedge \omega<0.$$
If $X$ admits a Hermitian metric $\omega$ with non-negative
Riemannian scalar curvature, then by formula (\ref{key2}), there
exists a Gauduchon metric with positive total Chern scalar
curvature. This is a contradiction.
 \eproof

\vskip 2\baselineskip

\section{Appendix: The scalar curvature relation on compact complex manifolds}
Let's recall some elementary settings (e.g. \cite[Section~2]{LY14}).
Let $(M, g, \nabla)$ be a $2n$-dimensional Riemannian manifold with
the Levi-Civita connection $\nabla$. The tangent bundle of $M$ is
also denoted by $T_\R M$. The Riemannian curvature tensor of
$(M,g,\nabla)$ is $$
R(X,Y,Z,W)=g\left(\nabla_X\nabla_YZ-\nabla_Y\nabla_XZ-\nabla_{[X,Y]}Z,W\right)$$
for tangent vectors $X,Y,Z,W\in T_\R M$. Let $T_\C M=T_\R M\ts \C$
be the complexification. We can extend the metric $g$ and the
Levi-Civita connection $\nabla$ to $T_{\C}M$ in the $\C$-linear way.
Hence for any $a,b,c,d\in \C$ and $X,Y,Z,W\in T_\C M$, we have  $$
R(aX,bY,cZ, dW)=abcd\cdot R(X,Y,Z,W).$$

\noindent Let $(M,g,J)$ be an almost Hermitian manifold, i.e.,
$J:T_\R M\>T_\R M$ with $J^2=-1$, and
 for any $X,Y\in T_\R M$, $ g(JX,JY)=g(X,Y)$. The Nijenhuis tensor $N_J:\Gamma(M,T_\R M)\times \Gamma(M,T_\R
M)\>\Gamma(M,T_\R M)$ is defined as $$
N_J(X,Y)=[X,Y]+J[JX,Y]+J[X,JY]-[JX,JY].$$ The almost complex
structure $J$ is called \emph{integrable} if $N_J\equiv 0$ and then
we call $(M,g,J)$ a Hermitian manifold. We can also extend $J$ to
$T_\C M$ in the $\C$-linear way. Hence for any $X,Y\in T_\C M$, we
still have $ g(JX,JY)=g(X,Y).$
 By Newlander-Nirenberg's
theorem, there exists a real coordinate system $\{x^i,x^I\}$ such
that $z^i=x^i+\sq x^I$ are  local holomorphic coordinates on $M$.
 Let's define a Hermitian form $h:T_\C M\times T_\C M\>\C$
by \beq h(X,Y):= g(X,Y),
 \ \ \ \ \ X,Y\in T_\C M.\label{complex}\eeq
By $J$-invariant property of $g$, \beq h_{ij}:=h\left(\frac{\p}{\p
z^i},\frac{\p}{\p z^j}\right)=0, \qtq{and} h_{\bar i\bar
j}:=h\left(\frac{\p}{\p \bar z^i},\frac{\p}{\p \bar
z^j}\right)=0\eeq and \beq h_{i\bar j}:=h\left(\frac{\p}{\p
z^i},\frac{\p}{\p \bar z^j}\right)=\frac{1}{2}\left(g_{ij}+\sq
g_{iJ}\right).\label{1003}\eeq It is obvious that $(h_{i\bar j})$ is
a positive Hermitian  matrix. Let $\omega$ be the fundamental
$2$-form associated to the $J$-invariant metric $g$:
 \beq \omega(X,Y)=g(JX,Y).\eeq
   In local complex coordinates,
 \beq \omega=\sq h_{i\bar j} dz^i\wedge d\bar z^j.\eeq

 \noindent
In the local holomorphic coordinates $\{z^1,\cdots, z^n\}$ on $M$,
the complexified Christoffel symbols are given by \beq
\Gamma_{AB}^C=\sum_{E}\frac{1}{2}g^{CE}\big(\frac{\p g_{AE}}{\p
z^B}+\frac{\p g_{BE}}{\p z^A}-\frac{\p g_{AB}}{\p
z^E}\big)=\sum_{E}\frac{1}{2}h^{CE}\big(\frac{\p h_{AE}}{\p
z^B}+\frac{\p h_{BE}}{\p z^A}-\frac{\p h_{AB}}{\p z^E}\big)
 \eeq
where $A,B,C,E\in \{1,\cdots,n,\bar{1},\cdots,\bar{n}\}$ and
$z^{A}=z^{i}$ if $A=i$, $z^{A}=\bar{z}^{i}$ if $A=\bar{i}$. For
example \beq \Gamma_{ij}^k=\frac{1}{2}h^{k\bar \ell}\left(\frac{\p
h_{j\bar \ell}}{\p z^i}+\frac{\p h_{i\bar \ell}}{\p z^j}\right),\
\Gamma_{\bar ij}^k=\frac{1}{2}h^{k\bar \ell}\left(\frac{\p h_{j\bar
\ell}}{\p \bar z^i}-\frac{\p h_{j\bar i}}{\p \bar
z^\ell}\right).\label{gamma} \eeq We also have $\Gamma_{\bar i\bar
j}^k=\Gamma_{ij}^{\bar k}=0$ by the Hermitian property
$h_{pq}=h_{\bar i\bar j}=0$.
 The complexified curvature components are \beq
R_{ABC}^D=\sum_ER_{ABCE}h^{ED}=-\left(\frac{\p\Gamma_{AC}^D}{\p
z^B}-\frac{\p\Gamma_{BC}^D}{\p
z^A}+\Gamma_{AC}^F\Gamma_{FB}^D-\Gamma_{BC}^F\Gamma_{AF}^D\right).
\eeq  By the Hermitian property again, we have \beq
 R_{i\bar jk}^{l}=-\left(\frac{\p
\Gamma^{l}_{ik}}{\p \bar z^j}-\frac{\p \Gamma^{l}_{\bar jk}}{\p
z^i}+\Gamma_{ ik}^{s}\Gamma^{l}_{\bar js}-\Gamma_{ \bar
jk}^{s}\Gamma^{l}_{ is}-{\Gamma_{\bar j k}^{\bar s}\Gamma_{i\bar
s}^l}\right).\label{cur}\eeq

 \vskip
1\baselineskip

\noindent It is computed in \cite[Lemma~7.1]{LY14} that \blemma On
the Hermitian manifold $(M,h)$, the Riemannian Ricci curvature of
the Riemannian manifold $(M,g)$ satisfies \beq Ric(X,Y)=h^{i\bar
\ell}\left[R\left(\frac{\p}{\p z^i}, X,Y, \frac{\p}{\p\bar
z^\ell}\right)+R\left(\frac{\p}{\p z^i}, Y,X, \frac{\p}{\p\bar
z^\ell}\right)\right]\label{ricci}\eeq for any $X,Y\in T_\R M$. The
Riemannian scalar curvature is \beq s=2h^{i\bar j}h^{k\bar
\ell}\left(2R_{i\bar \ell k\bar j}-R_{i\bar j k\bar
\ell}\right).\label{scalar}\eeq \elemma

\vskip 1\baselineskip \noindent The following result is established
in \cite[Corollary~4.2]{LY14} (see also some different versions in
\cite{Ga3}). For readers' convenience we include a straightforward
proof without using ``normal coordinates". \blemma\label{sc} On a
compact Hermitian manifold $(M,\omega)$, the Riemannian scalar
curvature $s$ and the Chern scalar curvature $s_{\mathrm C}$ are
related by \beq s=2s_{\mathrm
C}+\left(\la\p\p^*\omega+\bp\bp^*\omega,\omega\ra-2|\p^*\omega|^2\right)-\frac{1}{2}|T|^2,\eeq
where $T$ is the torsion tensor with
$$T_{ij}^k=h^{k\bar\ell}\left(\frac{\p  h_{j\bar \ell}}{\p z^i}-\frac{\p h_{i\bar\ell}}{\p z^j}\right).$$

 \elemma \bproof For simplicity, we
denote by
$$s_\text{R}=h^{i\bar j}h^{k\bar
\ell}R_{i\bar \ell k\bar j} \qtq{and} s_\text{H}=h^{i\bar j}h^{k\bar
\ell}R_{i\bar j k\bar \ell}.$$ Then, by formula (\ref{scalar}), we
have $s=4s_\text{R}-2s_\text{H}.$ In the following, we shall show
 \beq s_{\text{H}}=s_{\text{C}}-\frac{1}{2}\la
\p\p^*\omega+\bp\bp^*\omega,\omega\ra-\frac{1}{4}|T|^2\label{A2}\eeq
and \beq
s_\text{R}=s_\text{C}-\frac{1}{2}|\p^*\omega|^2-\frac{1}{4}|T|^2.\label{A1}\eeq

 It is easy to show that \beq \bp^*\omega=2\sq
\overline{\Gamma_{\bar{i}k}^k}dz^i \label{key11}\eeq and so
\beq-\frac{\p\p^*\omega+\bp\bp^*\omega}{2}=\sq\left(\frac{\p\Gamma_{\bar
j k}^k}{\p z^i}+\frac{\p\overline{\Gamma_{\bar{i}k}^k}}{\p\bar
z^j}\right)dz^i\wedge d\bar z^j.\eeq On the other hand, by formula
(\ref{cur}), we have \beq  {R}_{i\bar jk}^k=-\frac{\p
\Gamma^{k}_{ik}}{\p \bar z^j}+\frac{\p \Gamma^{k}_{\bar jk}}{\p
z^i}+{\Gamma_{\bar j k}^{\bar s}\Gamma_{i\bar s}^k}.\eeq A
straightforward calculation shows
$$h^{i\bar j} {\Gamma_{\bar j k}^{\bar s}\Gamma_{i\bar s}^k}=-\frac{1}{4}|T|^2.$$
 Moreover, we have \beq \left(-\frac{\p
\Gamma^{k}_{ik}}{\p \bar z^j}+\frac{\p \Gamma^{k}_{\bar jk}}{\p
z^i}\right)-\left(\frac{\p\Gamma_{\bar j k}^k}{\p
z^i}+\frac{\p\overline{\Gamma_{\bar{i}k}^k}}{\p\bar
z^j}\right)=-\frac{\p \Gamma^{k}_{ik}}{\p \bar
z^j}-\frac{\p\overline{\Gamma_{\bar{i}k}^k}}{\p\bar
z^j}=-\frac{\p^2\log \det(g)}{\p z^i\p\bar z^j}\eeq where the last
identity follows from (\ref{gamma}). Indeed, we have
$$\Gamma^{k}_{ik}+\overline{\Gamma_{\bar{i}k}^k}=h^{k\bar\ell}\frac{\p h_{k\bar\ell}}{\p z^i}=\frac{\p\log\det(g)}{\p z^i}.$$
 Hence, we obtain
$$s_\text{H}+\left\la \frac{\p\p^*\omega+\bp\bp^*\omega}{2}, \omega\right\ra=s_{\text{C}}-\frac{1}{4}|T|^2$$
which proves (\ref{A2}). Similarly, one can show (\ref{A1}).
 \eproof

\end{document}